\documentclass{article}

 \usepackage{graphicx}
 \usepackage[english]{babel}

 \usepackage{amsfonts,amssymb}
 \usepackage{amsmath}
 \usepackage{amsthm}
 \usepackage{cite}
 \usepackage[active]{srcltx}

 \newcommand{\tit}[1]{\begin{center}{\bf{\Large #1}}\end{center}}
 \newcommand{\aut}[1]{\centerline{{\bf #1}}}
 
 \newcommand{\email}[1]{\centerline{{\small e-mail: #1}}\vspace{\baselineskip}}

 \def\@evenhead{\vbox{\hbox to \textwidth{\thepage\hfil\sl\leftmark\strut}\hrule}}
 \def\@oddhead{\vbox{\hbox to \textwidth{\rightmark\hfill\thepage\strut}\hrule}}

\begin{document}
 \sloppy

\tit{Existence and Destruction of the Kantorovich Main Continuous Solutions of Nonlinear Integral Equations}

\aut{Denis Sidorov}

\section*{Introduction}

Consider the nonlinear Volterra integral equation of the second kind
\begin{equation}
x(t) = \int\limits_0^t K(t,s,x(s))\, ds, \,\, 0 < t < T<\infty,\, x(0)=0.  \label{sid1} 
\end{equation}

\noindent{\bf Definition}
[\cite{Kantor}, c.467]
Continuous function $x(t),$ satisfying the equation (\ref{sid1}), we name 
{\it Kantorovich main} solution if the sequence
 $$x_n = \int\limits_0^t K(t,s, x_{n-1}(s)) \, ds, \, \, x_0(t)=0$$
converges to function $x(t)$  $\forall t \in (0,T).$
If in addition $\lim\limits_{t \rightarrow T} |x(t)| = +\infty, $ then 
solution has blow-up point in the point $T.$

Let us find the guaranteed interval $\,[0,T)\,$
where exists the main solution such as the blow-up may occur
if one continue solution onto $[T, +\infty).$
Beside, one must find the positive continuous function $\hat{x}(t),$
defined on $[0,T)$
such as for the main solution $x(t)$
the following a priory estimate $|x(t)| \leq \hat{x}(t)$ if fulfilled for
 $t \in [0,T).$

In this paragraph we employ the classical approach by L.V.Kantorovich \cite{Kantor} 
(see chapter 12). For the majorizing equations construction we will use the algorithm
proposed in \cite{BanJSidSid}.

\noindent Let us introduce the following conditions\\
{\bf A.} Let function $K(t,s,x)$ be defined, continuous and differentiable wrt $x$
in $D = \{ 0<s<t<T, T\leq \infty, |x|<\infty \}.$\\
{\bf B.} Let we can construct functions
$m(s), \, \gamma(x)$ which are continuous, positive and monotonically increasing functions 
defined for $0<s<\infty, \, 0<x<\infty,$
such as in  $D$ for any $t$ from $[0, \infty)$
the following inequalities are fulfilled
$$|K(t,s,x)| \leq m(s) \gamma(|x|), $$
$$|K^{\prime}(t,s,x)| \leq m(s) \gamma^{\prime}(|x|). $$
The case of $\gamma(0)=0$ we exclude since in that case 
equation (\ref{sid1}) has only trivial solution.
Such solution is the main according to Kantorovich definition.
Below the functions $m(s), \, \gamma(x)$
are assumed positive, monotone increasing, and $\gamma(x)$
is assumed convex wrt  $x.$

\section{Integral Majorizing Equation}

Let us introduce majorizing integral equation
\begin{equation}
\hat{x}(t) = \int\limits_0^t m(s) \gamma(\hat{x}(s)) \, ds, \label{sid2} 
\end{equation}
 which is equivalent to the Cauchy problem for the differential equation with separable variables:
\begin{equation} \label{sid3}
    \left\{ \begin{array}{ll}
         \mbox{$\frac{d \hat{x}}{dt} = m(t) \gamma(\hat{x}(t))$} \\
         \mbox{$\hat{x}|_{t=0} = 0$}. \\
        \end{array} \right. 
\end{equation}
Thus, the solution $\hat{x}(t)$ of integral equation (\ref{sid2}) satisfies the equation $\Phi(x) = M(t),$
where $\Phi(x) = \int_{0}^x \frac{dx}{\gamma(x)},$
$M(t) = \int_0^t m(t) \, dt.$
Because of monotone increasing positive continuous function
$\Phi(x)$ exists inverse mapping$\Phi^{-1}$
with define area $[0,\infty),$ if $\lim\limits_{x\rightarrow \infty} \int\limits_0^x \frac{dx}{\gamma(x)} = +\infty,$
and with define area $[0, l), $ if $\lim\limits_{x\rightarrow \infty} \int\limits_0^x \frac{dx}{\gamma(x)} = l.$
Thus in the first case the Cauchy problem has unique positive solution
$\hat{x}(t)$ в $\mathcal{C}^{(1)}_{[0,\infty)},$
and in $\mathcal{C}^{(1)}_{[0,l)}$ in the second case.

\noindent {\bf Remark} \label{sidrem1}
Solution $\hat{x}(t)$ can be constructed with successive approximations
as solution of $\Phi(x)-M(t)=0.$
Indeed, equation $\Phi(t)-M(t)=0$
defines $\hat{x}(t)$ as explicit continuous function $\hat{x}(t) \rightarrow 0$  for  $t\rightarrow 0,$
since  $$\frac{d}{dx} \left( \Phi(x) - M(t)  \right) = \frac{1}{\gamma(x)} \neq 0. $$
Hence the solution $\hat{x}(t)$ can be constructed with successive
approximations:
$$x_n(t) = x_{n-1}(t) - \gamma(0) \left[ \Phi(x_{n-1}(t)) - M(t)  \right],\,
x_0=0$$
on the small interval $[0,\Delta], \, \Delta>0.$

Constructed solution can be continued on the whole domain of $\Phi^{-1}$ by repeated application of the implicit theorem application.

Let $\gamma(x)$ be polynomial with positive coefficients and $\gamma(0) \neq 0.$
In such case (ref. \cite{ilin}, p.344) function $\Phi(x)$ can be explicitly constructed
in terms of logarithms, arctangents and rational functions, which allows us in basic
cases to construct  $\Phi^{-1}$
and to explicitly build  $\hat{x}(t)$.

In general case in order to build $\hat{x}(t)=\Phi^{-1}(M(t)),$ satisfying equation
(\ref{sid2}), one may employ the following Lemma.

\noindent {\bf Lemma 1.}
\label{lem1}
If $\lim\limits_{x\rightarrow +\infty} \int\limits_0^x \frac{dx}{\gamma(x)}=+\infty,$
then majorizing equation (\ref{sid2}) has for $t \in [0, \infty)$
continuous solution $\hat{x}(t).$
Moreover, the sequence  
\begin{equation}
x_n(t) = \int\limits_0^t m(s) \gamma(x_{n-1}(s)) \, ds, \,\, x_0(t)=0
\label{sid4}
\end{equation}
converges for $t \in [0, \infty)$
to function $\hat{x}(t).$

\begin{proof}
Existence of the solution $\hat{x}(t)$ of equation (\ref{sid2})
follows from above proved solution existence of equivalent Cauchy problem (\ref{sid3}).
Herewith the sequence $\{x_n(t)\}$ will be monotone increasing
and  bounded above since
$\hat{x}(t)$ satisfies the equation (\ref{sid1}).
Hence  $\{x_n(t)\}$ has limit. Thus, the sequence
$\{x_n(t)\}$ is fundamental in the space
 $\mathcal{C}_{[0,T_1]}$ for
$T_1<\infty.$ The space  $\mathcal{C}_{[0,T_1]}$ is complete one and
convergence is uniform for $T_1< \infty$.

\end{proof}

\noindent {\bf Lemma 2.} \label{lem2}
Let $\lim\limits_{x\rightarrow +\infty} \int_0^x \frac{dx}{\gamma(x)} = l.$
Let us introduce the interval $[0,T_1),$ where $T_1>0$
in uniquely defined from the condition $\int_0^{T_1} m(s) \, ds = l.$
Then Cauchy problem (\ref{sid3})
has positive solution $\hat{x}(t) \in {\mathcal{C}}_{[0,T_1)},$
sequence $\{x_n(t)\}$ (см. (\ref{sid4}))
converges to $\hat{x}(t)$ as $n\rightarrow \infty,$ $\lim\limits_{t\rightarrow T_1}\hat{x}(t) = \infty.$

Proof of existence of the desired function $\hat{x}(t)$
on the interval $[0,T_1)$
follows from proved existence of inverse mapping  $\Phi^{-1}: [0,l) \rightarrow [0, +\infty)$. 
Since  $l=\int\limits_0^{\infty} \frac{dx}{\gamma(x)} = \int\limits_0^{T_1} m(s) \, ds,
$ based on Lemma 2 we have
$\lim\limits_{t\rightarrow T_1-0} \hat{x}(t) = +\infty.$

\noindent {\bf Theorem 1.}
\label{theosid1}
Let conditions {\bf (A)}, {\bf (B)} and $\lim\limits_{x\rightarrow +\infty} \int_0^x \frac{dx}{ \gamma(x)} = + \infty$ be fulfilled.
Then integral equation (\ref{sid1}) has main solution $x(t),$
defined for $t \in [0, \infty)$
and the following a priory estimate  $|x(t)|\leq \hat{x}(t)$ is valid, where 
 $\hat{x}(t)$ is solution of the Cauchy problem (\ref{sid3}).

\begin{proof}
Let us introduce the sequence
$$x_n(t) = \int\limits_0^t K(t,s, x_{n-1}(s))\, ds, \, x_0(t) \equiv 0, $$
$$\hat{x}_n(t) = \int\limits_0^t m(s) \gamma(\hat{x}_{n-1}(s))\, ds, \, \hat{x}_0(t) \equiv 0. $$
Then because of  \cite{BanJSidSid} and due to the Theorem conditions
we have the following inequality
$|x_{n+p}(t)-x_n(t)|\leq \hat{x}_{n+p}(t) - \hat{x}_n(t),$
$|x_n(t)| \leq \hat{x}_n(t)$
for $t \in [0,T_1], T_1 < \infty.$
Because of the Lemma 1 the positive monotonic increasing 
sequence $\hat{x}_n(t)$  is the Cauchy sequence
in the norm of the space ${\mathcal C}_{[0,T_1]},$
i.e. $\max\limits_{0\leq t \leq T_1} (x_{n+p}(t) - x_n(t)) \leq \varepsilon$
for $n \geq N(\varepsilon)$ and for arbitrary $p.$
Hence, $\max\limits_{0\leq t \leq T_1} |x_{n+p}(t) - x_n(t)| \leq \varepsilon$
for $n\geq N(\varepsilon), \, \forall p.$ Hence the sequence $\{x_n(t)\}$
is in the sphere $S(0,\hat{x}(t))$ and its the Cauchy sequence . 
Because of completeness of the space  $\mathcal{C}_{[0,T_1]}$\\ $\exists \lim\limits_{n\rightarrow \infty} x_n(t) = x(t).$  
And $|x(t)| \leq \hat{x}(t).$
Since $K(t,s,x)$ is continuous wrt $x$ then function $x(t)$ satisfies the condition (\ref{sid1}). 
The theorem is proved.

The classic Hartman-Wintner theorem  on the Cauchy problem solution continuation on semi-axis \cite{Yudovich} follows from
the proved theorem. 
\end{proof}

\noindent{\bf Corollary 1}
Let $\gamma(x) = a +bx$ in the condition (B) be linear function, $a \geq 0, b>0.$
Then for main solution of equation (\ref{sid1}) the following a priory estimate $$|x(t)|\leq 
 {a} \int\limits_0^t 
m(z) \exp   \left(b \int\limits_z^t m(s) \, ds \right)  \, dz $$
is fulfilled for $0\leq t <\infty.$

In order to prove this corollary it is enough to verify
that the Cauchy problem
$\frac{dx}{dt} = m(t)(a+bx(t)), \, x|_{t=0}=0$
has the solution  $$\hat{x}(t) = {a} \int\limits_0^t 
m(z) \exp   \left(b \int\limits_z^t m(s) \, ds \right)  \, dz$$
for $0\leq t <\infty.$

\noindent {\bf Theorem 2.} \label{theo2}
Let conditions (A) and (B) be fulfilled.
Let $\lim\limits_{x\rightarrow \infty} \int\limits_0^x \frac{dx}{\gamma(x)} = l.$
Introduce the interval  $[0,T_1],$  where $ T_1>0$ is uniquely defined 
from the equality $\int\limits_0^{T_1} m(s) \, ds =l.$
Then the integral equation (\ref{sid1}) in  ${\mathcal{C}}_{[0,T_1)}$ has
the main solution $x(t).$ For  $t\in [0,T_1)$ its a priory estimate $|x(t)|\leq \hat{x}(t),$
is fulfilled, where
 $\hat{x}(t)$ is the solution of the Cauchy problem  (\ref{sid3}), $\lim\limits_{t\rightarrow T_1} \hat{x}(t) = +\infty.$

Proof follows from the Theorem 1 taking into account 
the Lemma 2 results. 

\noindent{\bf Corollary 2. }
(Alternative global solvability of equation  (\ref{sid1}))
Let conditions of Theorem 2 be fulfilled. Then either solution of the equation (\ref{sid1})
can be continued on the whole semi-axis or on  $[T_1, +\infty)$ or there the 
blow-up point exists.
 
Corollary 2 refines and generalizes the V.~I.~Yudovich theorem on
global resolvability of the Cauchy problem  (ref. Theorem 1 in \cite{Yudovich}, p. 19).\\

\noindent{\bf Corollary 3. }
Let conditions of the Theorem 2 be fulfilled, and in addition let $\gamma(x) = ax^2+bx+c,$ 
where $4ac-b^2>0,$ $a >0, b >0, c \geq 0.$
Then main solution $x(t)$ of the equation (\ref{sid1}) exists on  $[0,T_1],$
where positive $T_1$ is defined from condition
$$\int\limits_0^{T_1} m(s) \, ds = \frac{2}{\sqrt{4ac-b^2}} \left(\frac{\pi}{2}-arctg \frac{b}{\sqrt{4ac-b^2}}  \right). $$
For $t \in [0,T_1)$
we have the estimate $|x(t)| \leq \hat{x}(t),$
where function $\hat{x}(t)$ is defined by formula
\begin{equation}
\hat{x}(t) = \frac{\sqrt{4ac-b^2}}{2a} tg \biggl[ arctg \frac{b}{\sqrt{4ac-b^2}} + 
\label{sid14}
\end{equation}
$$
+\frac{\sqrt{4ac-b^2}}{2} \int\limits_0^t m(s)\, ds  \biggr ] -
\frac{b}{2a}.
$$
\begin{proof}

Solution to the corresponding Cauchy problem (\ref{sid3}) 
can be constructed easily when Corollary 3 conditions are fulfilled.
Indeed, (see \cite{Dvait}, p.36), $\Phi(x) = \int\limits_0^x \frac{dx}{ax^2 +bx +c} = \frac{2}{\sqrt{4ac-b^2}} 
\left ( arctg \frac{2ax+b}{\sqrt{4ac-b^2}} - arctg \frac{b}{\sqrt{4ac-b^2}} \right ),$
$M(t)= \int_0^t m(s)\, ds.$
As results the equality $\int\limits_0^{\infty} \frac{dx}{\gamma(x)} = \int\limits_0^{T_1} m(s) \, ds,$
which serves to find $T_1,$
means that 
\begin{equation}
\frac{2}{\sqrt{4ac-b^2}} 
\left ( \frac{\pi}{2} - arctg \frac{b}{\sqrt{4ac-b^2}} \right ) = \int_0^{T_1} m(s) \, ds. 
\label{sid15}
\end{equation}
From equation $\Phi(x) = M(t)$ it follows that majorizing function $\hat{x}(t),$
on $[0,T_1),$ must be constructed with formula
(\ref{sid14}).
\end{proof}
Corollary 3 can be employed can be used in the problem of extending the solution of the equation (\ref{sid1}) with parameters. 

Indeed, let conditions of the Theorem 2 be fulfilled and function $K$ 
depends on parameter $\lambda \in \mathbb{R}^n,\, ||\lambda||<\delta$ (i.e. $K  = K(t,s,x,\lambda)$).
Then in condition {\bf B} function $\gamma(|x|, || \lambda||)$
will depend on this parameter's norm. Let  $\gamma(|x|,0)=0.$
Then equation (\ref{sid1}) has trivial solution for  $\lambda=0.$ 
If  $\lambda \neq 0$ equation  (\ref{sid1}) can has nontrivial main solution.
Next result allows us to estimate the interval where exists main
solution to the equation (\ref{sid1}) for $0< ||\lambda||<\delta.$\\

\noindent {\bf Corollary 4\, }
Let conditions of the Theorem 2 are fulfilled and let $m(s)=1,$  $\gamma(x,\lambda) = a(||\lambda||)x^2 +  b(||\lambda||)x+
c(||\lambda||).$
Let $a(||\lambda||),$  $b(||\lambda||),$  and $c(||\lambda||)$ are positive
infinitesimal functions for $||\lambda|| \rightarrow 0.$
Suppose that in a punctured neighborhood $0<||\lambda||<\delta$
the following inequalities are fulfilled
$4a(||\lambda||) c(||\lambda||) - b^2 (||\lambda||) = \Delta (||\lambda||)>0,$
$$ \sup\limits_{||\lambda||<\delta} \arctan \frac{b(||\lambda||)}{\sqrt{\Delta(||\lambda||)}} 
= \sigma, \, \sigma<\frac{\pi}{2}.$$

Then for $0<||\lambda||<\delta$ equation
(\ref{sid1}) has the main solution $x(t,\lambda),$
defined on the intervals $[0,T(||\lambda||),$
  $T(||\lambda||) = \frac{2}{\sqrt{\Delta(||\lambda||)}} \left(\frac{\pi}{2}-\sigma \right).$
Therefore $\lim\limits_{||\lambda||\rightarrow 0} T(||\lambda||)=+\infty$
and the following a propry estimate for the main solution is fulfilled
$|x(t,\lambda)| \leq \hat{x}(t,\lambda).$
Here the positive majorizing function $\hat{x}(t,\lambda)$
is defined on the interval $[0,T(||\lambda||)$
using the formula (\ref{sid14})
where coefficients are defined as follows: $a=a(||\lambda||), b=b(||\lambda||),$ $ c=c(||\lambda||),$
$\lim\limits_{t\rightarrow T(||\lambda||)} \hat{x}(t,\lambda)=+\infty.$


\section{Algebraic majorants}

In order to estimate the gauranteed closed interval $[0,T]$ for existence of main solution of the equation (\ref{sid1})
and its norm estimation ${\mathcal C}_{[0,T]}$ the algebraic majorants are useful.

Indeed, let condition (A) be fulfilled and let in addition the following condition be fulfilled:\\
$\,$\\
{\bf (С)} Let there exists continuous, differentiable and convex wrt $r,$  positive and monotonic
increasing function $M(\rho, s, r),$ defined for $\rho\geq 0, s \geq 0, r>0$
such as in the area $D$ for $t \in [0, \rho], \, |x| \leq r$ the following inequalities are fulfilled
$$\left|\int\limits_0^t K(t,s,x(s))\, ds\right| \leq \int\limits_0^{\rho} M(\rho,s,r) \, ds,$$
$$\left|\int\limits_0^t K^{\prime}(t,s,x(s))\, ds\right| \leq \int\limits_0^{\rho} M_{r}^{\prime}(\rho,s,r)\, ds.$$
Let us introduce function $$M(\rho, r) = \int\limits_0^{\rho} M(\rho,s,r)\, ds$$
and it's positive derivative $$ M_r^{\prime}(\rho, r) = \int\limits_0^{\rho} M_r^{\prime}(\rho,s,r)\, ds.$$

\noindent {\bf Lemma 3.} \label{lem3}
   System
\begin{equation}
    \left\{ \begin{array}{ll}
         \mbox{$r = M(r,\rho)$,} \\
         \mbox{$1=M_r^{\prime}(r,\rho)$} \\
        \end{array} \right.
\label{sid16} 
\end{equation}
has unique positive solution $r^*, \rho^*.$
Moreover, for any $\rho \in [0,\rho^*]$
the equation $r=M(r,\rho)$
has main solution $r(\rho),$ i.e. monotonic increasing sequence
$r_n = M(r_{n-1},\rho), r_0=0,$ converges to the solution of equation $r=M(r,\rho)$
for any $\rho \in [0,\rho^*].$

Proof is geometrically obvious (here readers may refer to \cite{Greb}, p. 218),
if on the plane $(y,r)$ one consider the graphs of the curves 
 $y=M(r,\rho)$ for various  $\rho$  and bisection $y=r.$
Line $y=r$ tangents the curve $y=M(r,\rho^*)$ in the point $(r^*, \rho^*),$

\noindent {\bf Theorem 3} \label{sidtheo3}
Let conditions (A) and (C) be fulfilled,  
 $(r^*, \rho^*)$ is positive solution of the system (\ref{sid16}).
Then main solution $x(t)$ of the equation (\ref{sid1}) exists in ${\mathcal{C}}_{[0,\rho^*]}$
and the following estimate is fulfilled $$\max\limits_{0\leq t\leq\rho^*}|x(t)| \leq r^*.$$

\begin{proof}
Let us introduce two sequences $$x_n(t) = \int\limits_0^t K(t,s,x_{n-1}(s)) \, ds,$$ 
$$r_n = M(r_{n-1},\rho^{\ast}), $$ где $ x_0(t)=0, \, r_0=0.$
Then the following inequality $||x_{n+p}-x_n||_{{\mathcal{C}}_{[0,\rho^*]}} \leq r_{n+p}-r_n$ is valid for
 $n \geq N(\varepsilon)$ and for any $p.$
Therefore since due to the Lemma 3 we have $\lim\limits_{n\rightarrow 0} r_n=r^*,$
then $r_{n+p}-r_n \leq \varepsilon$ for $n \geq N(\varepsilon)$ и для любых $p.$
Therefore the sequence $\{x_n(t)\},$  for $ x_0=0$ remains fundamental, $||x_n||\leq r^*$
and the theorem is proved.
\end{proof}

Let us consider the following example:
$$x(t) = \int\limits_0^t \left(K_2(t,s)x^2(s)+K_1(t,s)x(s)+K_0(t,s)   \right) s^2  \, ds. $$
Let $$\sup\limits_{0\leq s \leq t <\infty, \, i=1,2,3}  |K_i(t,s)|\leq 1.$$
The corresponding majorant algebraic system
$$
\left\{ \begin{array}{ll}
         \mbox{$r = \frac{\rho^3}{3} (1+r+r^2)$}, \\
         \mbox{$1 = \frac{\rho^3}{3}(1+2r)$} \\
        \end{array} \right. 
$$
has the following solution:  $r^*=1, \, \rho^*=1.$
Therefore, based on Theorem 3 the integral equation  has the main solution
 $x(t) \in \mathcal{C}_{[0,1]},$ $||x||\leq 1.$ 
The integral majorant 
$x(t) = \int_0^t (x(s)^2+x(s)) \,ds +\frac{t^3}{3}$ and Corollary 3 gives us
more complete information regarding the solution.
Indeed, the integral equation has continuous solution $x(t)$ on the interval $\left[0, 1.5365\right)$
and the following estimate is fulfilled
$|x(t)| \leq 0.8660 \tan \left (  0.5236 + 0.2886 t^3 \right ) - 0.5$ for $0\leq t < 1.5365.$

\begin{figure}[htbp]
\centering{\includegraphics[scale=0.6]{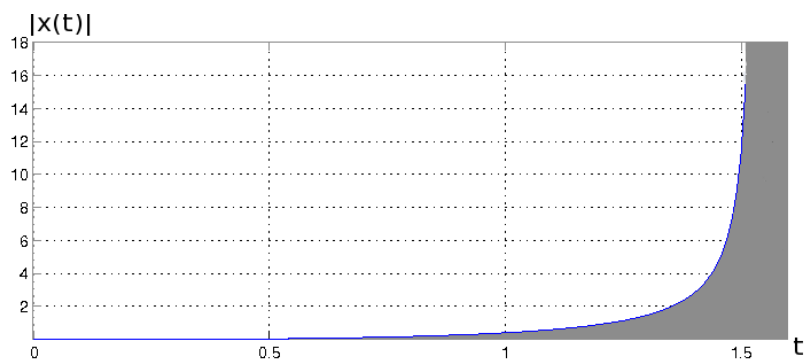}}
\caption{Majorant\, of the solution.\,\, Solution $|x(t)|$ belongs to the gray zone}\label{sidfig1}
\end{figure}

\end{document}